\documentclass[twocolumn]{autart}    

\usepackage{graphicx}          
\usepackage{graphics} 
\usepackage{epsfig} 
\usepackage{amsmath, amssymb, amsfonts}
\usepackage{natbib}
\usepackage{graphicx,color}

\newcommand\oprocendsymbol{\hbox{$\square$}}
\newcommand\oprocend{\relax\ifmmode\else\unskip\hfill\fi\oprocendsymbol}

\newtheorem{defi}{\bf{\emph{Definition}}}
\newtheorem{ex}{\bf{\emph{Example}}}
\newtheorem{asmp}{Assumption}

\def\E{{\mathbb E}}
\def\P{{\tilde P}}
\def\R{{\mathbb R}}
\def\N{{\mathcal N}}
\def\I{{\mathcal I}}
\def\z{{\tilde z}}
\def\u{{\tilde u}}
\def\y{{\tilde y}}
\def\x{{\tilde x}}
\def\e{{\tilde e}}
\def\tr{{\text{tr}}}
\def\det{{\text{det}}}

\renewcommand{\theenumi}{(\roman{enumi})}

\begin{document}

\begin{frontmatter}

  \title{Data-Injection Attacks in Stochastic Control Systems:
    Detectability and Performance Tradeoffs\thanksref{footnoteinfo}}


  \thanks[footnoteinfo]{Work supported in part by awards NSF
    ECCS-1405330 and ONR N00014-14-1-0816. Corresponding author:
    V. Gupta. Tel. +1 574 631 2294. A preliminary version of this work
    appeared in \cite{Bai15}. With respect to \cite{Bai15}, this paper
    (i) considers more general systems with multiple inputs and
    outputs, (ii) completes and extends technical proofs, and (iii)
    provides further insight into the design of optimal stealthy
    attacks in stochastic cyber-physical systems.}

\author[Vijay]{Cheng-Zong Bai}\ead{cbai@nd.edu},    
\author[Fabio]{Fabio Pasqualetti}\ead{fabiopas@engr.ucr.edu},               
\author[Vijay]{Vijay Gupta}\ead{vgupta2@nd.edu}  

\address[Vijay]{Department of Electrical Engineering, University of Notre Dame, Notre Dame, IN USA}        
\address[Fabio]{Department of
    Mechanical Engineering, University of California, Riverside, CA USA}             
          

\begin{abstract}                          
  Consider a stochastic process being controlled across a
  communication channel. The control signal that is transmitted across
  the control channel can be replaced by a malicious attacker. The
  controller is allowed to implement any arbitrary detection algorithm
  to detect if an attacker is present. This work characterizes some
  fundamental limitations of when such an attack can be detected, and
  quantifies the performance degradation that an attacker that seeks
  to be undetected or stealthy can introduce.
\end{abstract}

\end{frontmatter}
	
\section{Introduction}
Using communication channels to inject malicious data that degrades
the performance of a cyber-physical system has now been demonstrated
both theoretically and practically
\cite{JPF-RR:11,SK:03,Mo14,FP-FD-FB:10y,GR:08,JS-MM:07}.
Intuitively, there is a tradeoff between the performance degradation
an attacker can induce and how easy it is to detect the attack
\cite{Teixeira12}. Quantifying this tradeoff is of great interest to
operate and design secure cyber-physical systems (CPS).

As explained in more detail later, for noiseless systems, zero
dynamics provide a fundamental notion of stealthiness of an attacker,
which characterizes the ability of an attacker to stay undetected even
if the controller can perform arbitrary tests on the data it
receives. However, similar notions for stochastic systems have been
lacking. In this work, we consider stochastic cyber-physical systems,
propose a graded stealthiness notion, and characterize the performance
degradation that an attacker with a given level of stealthiness can
induce. The proposed notion is fundamental in the sense that we do not
constraint the detection test that the controller can employ to detect
the presence of an attack.

\noindent
\textbf{Related work} Security of cyber-physical systems is a growing
research area. Classic works in this area focus on the detection of
sensor and actuator failures in control systems \cite{RP-PF-RC:89},
whereas more recent approaches consider the possibility of intentional
attacks at different system layers; e.g., see \cite{FP-FD-FB:15}. Both
simple attacks, such as jamming of communication channels
\cite{Foroush13}, and more sophisticated attacks, such as replay and
data injection attacks, have been considered \cite{ym-bs:10a,rs:11}.

One way to organize the literature in this area is based on the
properties of the considered cyber-physical systems. While initial
studies focused on static systems
\cite{DG-HS:10,AG-EB-MG-MM-PK-KP:11,YL-MKR-PN:09,ahmr-alg:11,AT-AS-HS-KHJ-SSS:10},
later works exploited the dynamics of the system either to design
attacks or to improve the performance of the detector that a
controller can employ to detect if an attack is present
\cite{SB-TB:13,FH-PT-SD:11,SM-QZ-YZ-SG-TB:13,MM-QZ-TA-TB-JPH:11,MZ-SM:11,QZ:HT:TB:13}. For
noiseless cyber-physical systems, the concept of stealthiness of an
attack is closely related to the control-theoretic notion of zero
dynamics \cite[Section 4]{GB-GM:91}. In particular, an attack is
undetectable in noiseless systems if and only if it excites only the
zero dynamics of an appropriately defined input-output system
describing the system dynamics, the measurements available to the
security monitor, and the variables compromised by the attacker
\cite{Fawzi14,FP-FD-FB:10y}. 
For cyber-physical systems driven by noise, instead, the presence of
process and measurements noise offers the attacker an additional
possibility to tamper with sensor measurements and control inputs
within acceptable uncertainty levels, thereby making the detection
task more difficult.

Detectability of attacks in stochastic systems remains an open
problem. Most works in this area consider detectability of attacks
with respect to specific detection schemes employed by the controller, such as the classic bad
data detection algorithm \cite{Cui12,ym-bs:10a}. The trade-off between
stealthiness and performance degradation induced by an attacker has
also been characterized only for specific systems and detection
mechanisms \cite{Kosut11,Kwon13,Liu11,Mo14}, and a thorough analysis of
 resilience of stochastic control systems to arbitrary attacks is
still missing. While convenient for analysis, the restriction to a
specific class of  detectors prevents the characterization
of fundamental detection limitations. In our previous work
\cite{Bai14}, we proposed the notion of $\epsilon$-marginal
stealthiness to quantify the stealthiness level in an estimation
problem with respect to the class of ergodic detectors. In this work,
we remove the assumption of ergodicity and introduce a notion of stealthiness for stochastic control systems
that is independent of the attack detection algorithm, and thus
provides a fundamental measure of the stealthiness of attacks in
stochastic control systems. Further, we also characterize the
performance degradation that such a stealthy attack can induce.

We limit our analysis to linear, time-invariant plants with a
controller based on the output of an asymptotic Kalman filter, and to
injection attacks against the actuation channel only. Our choice of
using controllers based on Kalman filters is not restrictive. In fact,
while this is typically the case in practice, our results and analysis
are valid for arbitrary control schemes. Our choice of focusing on
attacks against the actuation channel only, instead, is motivated by
two main reasons. First, actuation and measurements channels are
equally likely to be compromised, especially in networked control
systems where communication between sensors, actuators, plant, and
controller takes place over wireless channels. Second, this case has
received considerably less attention in the literature -- perhaps due
to its enhanced difficulty -- where most works focus on attacks
against the measurement channel only; e.g., see
\cite{Fawzi14,AT-AS-HS-KHJ-SSS:10}. We remark also that our framework
can be extended to the case of attacks against the measurement
channel, as we show in \cite{Bai14} for scalar systems and a different
notion of stealthiness.
  
Finally, we remark that since the submission of this work, some recent
literature has appeared that builds on it and uses a notion of attack
detectability that is similar to what we propose in \cite{Bai14,Bai15}
and in this paper. For instance, \cite{EK-SD-LS:16} extends the notion
of $\epsilon$-stealthiness of \cite{Bai15} to higher order systems,
and shows how the performance of the attacker may differ in the scalar
and vector cases (in this paper we further extend the setup in
\cite{EK-SD-LS:16} by leveraging the notion of right-invertibility of
a system to consider input and output matrices of arbitrary
dimensions). In \cite{RZ-PV:16}, the authors extend the setup in
\cite{Bai15} to vector and not necessarily stationary systems, but
consider a finite horizon problem. In \cite{ZG-DS-KHJ-LS:16}, the
degradation of remote state estimation is studied, for the case of an
attacker that compromises the system measurements based on a linear
strategy. Two other relevant recent works are \cite{sw-bs-sk-ad:16}
that uses the notion of Kullback-Liebler divergence as a causal
measure of information flow to quantify the effect of attacks on the
system output, while \cite{YC-SK-JMFM:16} characterizes optimal attack
strategies with respect to a linear quadratic cost that combines
attackers control and undetectability goals.

\noindent
\textbf{Contributions} The main contributions of this paper are
threefold.  First, we propose a notion of $\epsilon$-stealthiness to
quantify detectability of attacks in stochastic cyber-physical
systems. Our metric is motivated by the Chernoff-Stein Lemma in
detection and information theories and is universal because it is
independent of any specific detection mechanism employed by the
controller. Second, we provide an information theoretic bound for the
degradation of the minimum-mean-square estimation error caused by an
$\epsilon$-stealthy attack as a function of the system parameters,
noise statistics, and information available to the attacker. Third, we
characterize optimal stealthy attacks, which achieve the maximal
degradation of the estimation error covariance for a stealthy
attack. For right-invertible systems \cite[Section 4.3.2]{GB-GM:91},
we provide a closed-form expression of optimal $\epsilon$-stealthy
attacks. The case of single-input single-output systems considered in
our conference paper~\cite{Bai15} is a special case of this
analysis. For systems that are not right-invertible, we propose a
sub-optimal $\epsilon$-stealthy attack with an analytical expression
for the induced degradation of the system performance. We include a
numerical study showing the effectiveness of our bounds. Our results
provide a quantitative analysis of the trade-off between performance
degradation that an attacker can induce versus a fundamental limit of
the detectability of the attack.

\noindent
\textbf{Paper organization} Section \ref{sec:sys_model} contains the
mathematical formulation of the problems considered in this paper. In
Section \ref{sec:stealthy}, we propose a metric to quantify the
stealthiness level of an attacker, and we characterize how this metric
relates to the information theoretic notion of Kullback-Leibler
Divergence. Section \ref{sec:fundamental_limits_given_W} contains the
main results of this paper, including a characterization of the
largest performance degradation caused by an $\epsilon$-stealthy
attack, a closed-form expression of optimal $\epsilon$-stealthy
attacks for right invertible systems, and a suboptimal class of
attacks for not right-invertible systems. Section \ref{sec:numerical}
presents illustrative examples and numerical results. Finally, Section
\ref{sec:conclusion} concludes the paper.

\section{Problem Formulation}\label{sec:sys_model}
\textbf{Notation:} The sequence $\{x_n\}_{n=i}^j$ is denoted by
$x_i^j$ (when clear from the context, the notation $x_i^j$ may also
denote the corresponding vector obtained by stacking the appropriate
entries in the sequence). This notation allows us to denote the
probability density function of a stochastic sequence $x_i^j$
$f_{x_i^j}$, and to define its differential entropy $h(x_i^j)$ as
\cite[Section 8.1]{Cover06}
\begin{align*}
  h(x_i^j)\triangleq\int_{-\infty}^\infty -f_{\tilde{x}_i^j}(t_i^j) \log
  f_{\tilde{x}_i^j}(t_i^j) d t_i^j .
\end{align*}
Let $x_1^k$ and $y_1^k$ be two random sequences with
probability density functions (pdf) $f_{x_1^k}$ and $f_{y_1^k}$,
respectively. The Kullback-Leibler Divergence (KLD)~\cite[Section 8.5]{Cover06} between $x_1^k$ and $y_1^k$ is defined as
\begin{equation}\label{eq:KLD}
  D\big( x_1^k \big\| y_1^k \big) \triangleq \int_{-\infty}^\infty \log \frac{f_{x_1^k}(t_1^k)}{f_{y_1^k}(t_1^k)} f_{x_1^k}(t_1^k) d t_1^k.
\end{equation}
The KLD is a non-negative quantity that gauges the dissimilarity
between two probability density functions with
$D\big( x_1^k \big\| y_1^k \big)=0$ if $f_{x_1^k}=f_{y_1^k}$. Also,
the KLD is generally not symmetric, that is,
$D\big( x_1^k \big\| y_1^k \big) \neq D\big( y_1^k \big\| x_1^k
\big)$.
A Gaussian random vector $x$ with mean $\mu_x$ and covariance matrix
$\Sigma_x$ is denoted by $x\sim\N(\mu_x, \Sigma_x)$. We let $I$ and
$O$ be the identity and zero matrices, respectively, with their
dimensions clear from the context. We also let
$\mathbb{S}^{n}_{+}$ and $\mathbb{S}^{n}_{++}$ denote the sets of
$n\times n$ positive semidefinite and positive definite matrices,
respectively. For a square matrix $M$, $\tr(M)$ and $\det(M)$ denote
the trace and the determinant of $M,$ respectively.

\begin{figure}
    \centering
    \includegraphics[width=1\columnwidth]{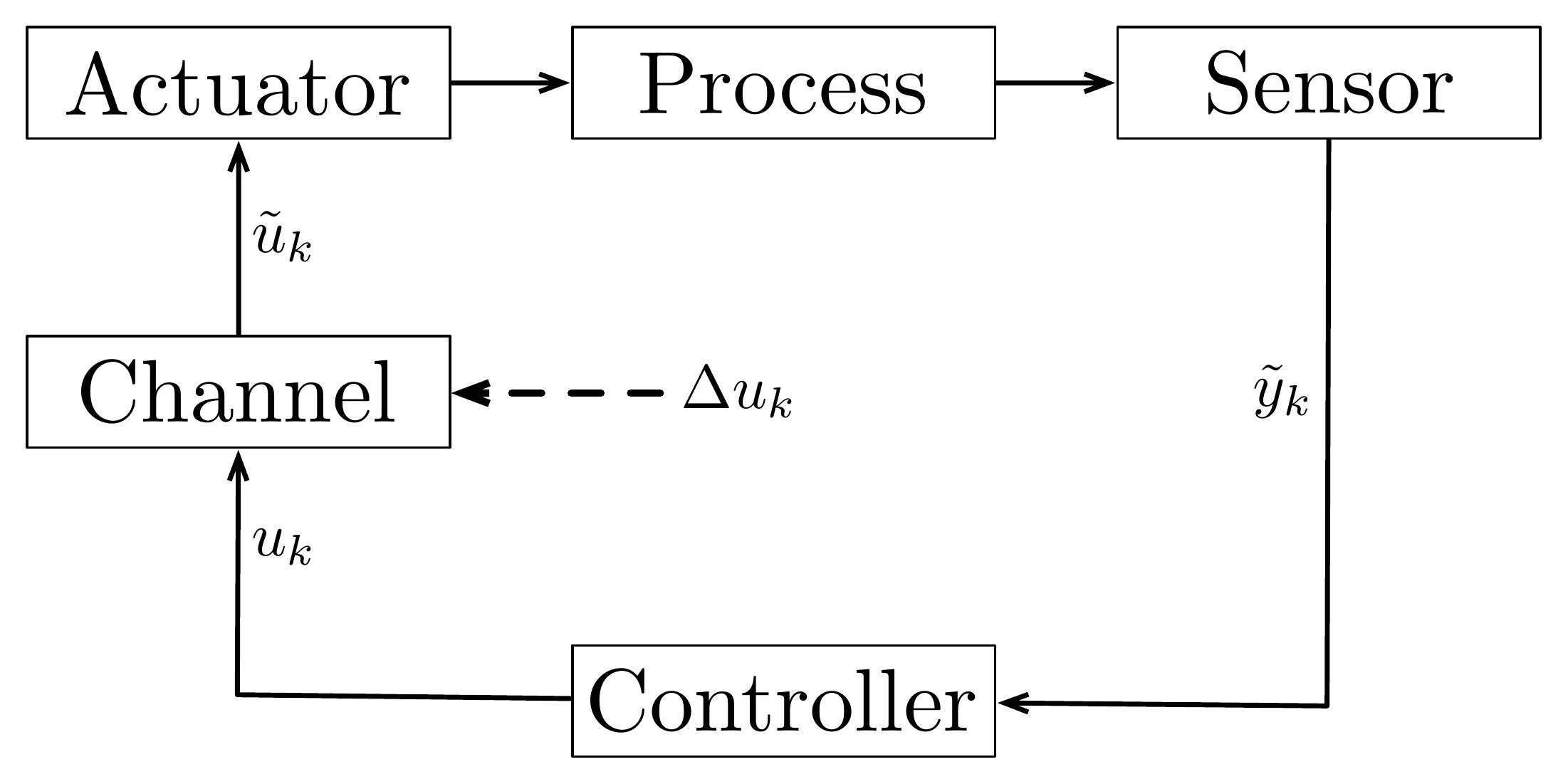}
    \caption{Problem setup considered in the paper.}
    \label{fig: setup}
\end{figure}

We consider the setup shown in Figure~\ref{fig: setup} with the
following assumptions:

\textbf{Process:} The process is described by the following linear
time-invariant (LTI) state-space representation:
\begin{equation}\label{eq:LTI}
  \begin{aligned}
    x_{k+1} &= A x_k + B u_k + w_k,\\
    y_k &= C x_k + v_k,
  \end{aligned}
\end{equation}
where $x_{k}\in\R^{N_x}$ is the process state, $u_{k}\in\R^{N_{u}}$ is
the control input, $y_{k}\in\R^{N_{y}}$ is the output measured by the
sensor, and the sequences $w_1^\infty$ and $v_1^\infty$ represent
process and measurement noises, respectively. 
%
\begin{asmp}
  The noise random processes are independent and identically
  distributed (i.i.d.)  sequences of Gaussian random vectors with
  $w_k\sim\N(0,\Sigma_w)$, $v_k\sim\N(0, \Sigma_v)$, 
  $\Sigma_w \in \mathbb{S}^{N_x}_{++}$, and
  $\Sigma_v \in \mathbb{S}^{N_y}_{++}$. 
\end{asmp}

\begin{asmp}
\label{asmp_zero}
The state-space realization $(A, B, C)$ has no invariant zeros
\cite[Section 4.4]{GB-GM:91}. In particular, this assumption implies
that the system $(A,B,C)$ is both controllable and
observable. 
\end{asmp}
\begin{asmp}\label{asmp_kalman} 
  The controller uses a Kalman filter to estimate and monitor the
  process state. Note that the control input itself may be calculated
  using an arbitrary control law. The Kalman filter, which calculates
  the Minimum-Mean-Squared-Error (MMSE) estimate $\hat{x}_{k}$ of
  $x_{k}$ from the measurements $y_1^{k-1}$, is described as
\begin{equation}\label{eq:kf}
  \hat x_{k+1} = A \hat{x}_k + K_k(y_k - C \hat{x}_{k}) + B u_k ,
\end{equation}
where the Kalman gain $K_k$ and the error covariance matrix
$P_{k+1}\triangleq
\E\big[(\hat{x}_{k+1}-x_{k+1})(\hat{x}_{k+1}-x_{k+1})^T\big]$
are calculated through the recursions
\begin{align*}
  K_k &= A P_k C^T (C P_K C^T + \Sigma_v)^{-1} , \text{ and}\\
  P_{k+1} &= A P_k A^T - A P_k C^T (C P_k C^T + \Sigma_v)^{-1}C
            P_k A^T  + \Sigma_w ,
\end{align*}
with initial conditions $\hat{x}_1 = \E[x_1]=0$ and
$P_1=\E[x_1 x_1^T]$.  
\end{asmp}
\begin{asmp}
  Given Assumption~\ref{asmp_zero}, $\lim_{k\to\infty}P_k=P$, where
  $P$ is the unique solution of a discrete-time algebraic Riccati
  equation. For ease of presentation, we assume that $P_1=P$, although
  the results can be generalized to the general case at the expense of
  more involved notation. Accordingly, we drop the time index and
  let $K_k=K$ and $P_k=P$ at every time step $k$. Notice that this
  assumption also implies that the innovation sequence $z_1^\infty$
  calculated as $z_k\triangleq y_k-C\hat{x}_k$ is an i.i.d. Gaussian
  process with $z_k\sim \N(0, \Sigma_z)$, where
  \begin{equation}\label{eq:Sigma_z}
    \Sigma_z = CPC^T + \Sigma_v\in\mathbb{S}_{++}^{N_y}.
  \end{equation}
\end{asmp}

Let $G(\mathcal{Z})$ denote the $N_y\times N_u$ matrix transfer
function of the system $(A,B,C)$. We say that the system $(A,B,C)$ is
right invertible if there exists an $N_u\times N_y$ matrix transfer
function $G_{RI}(\mathcal{Z})$ such that
$G(\mathcal{Z})G_{RI}(\mathcal{Z}) = I_{N_y}$.  
	
\textbf{Attack model:} An attacker can replace the input sequence
$u_1^\infty$ with an arbitrary sequence $\u_1^\infty$. Thus, in the
presence of an attack, the system dynamics are given by
\begin{align}\label{eq:LTI_attack}
  \tilde x_{k+1} &= A \tilde x_k + B\u_k + w_k, \notag\\
  \y_k &= C \tilde x_k + v_k.
\end{align}        
Notice that the sequence of measurements $\y_1^\infty$ generated by
the sensor in the presence of an attack $\u_1^\infty$ is different
from the nominal measurement sequence $y_1^\infty$.  We assume that
the attacker knows the system parameters, including the matrices $A$,
$B$, $C$, $\Sigma_w$, and $\Sigma_v$. The attack input $\u_1^\infty$
is constructed based on the system parameters and the
\emph{information pattern} $\I_k$ of the attacker. We make the following assumptions on the attacker's information pattern:
\begin{asmp}\label{assm_A1} 
  The attacker knows the control input $u_k$; thus $u_k \in \I_k$ at
  all times $k$. Additionally, the attacker does not know the noise
  vectors for any time. 
\end{asmp}
\begin{asmp}\label{assm_A2}	
  The attacker has perfect memory; thus, $\I_k\subseteq\I_{k+1}$ at
  all times $k$. 
\end{asmp}
\begin{asmp}\label{assm-A3}
  The attacker has causal information; in particular, $\I_k$ is
  independent of $w_k^\infty$ and $v_{k+1}^\infty$ for all
  $k$. 
\end{asmp}

\begin{ex}{\bf \emph{(Attack scenarios)}}
Attack scenarios satisfying
Assumptions~\ref{assm_A1}-\ref{assm-A3}  include the cases when:
\begin{enumerate}
\item the attacker knows the control input exactly, that is,
  $\I_k=\{u_1^k\}$.
\item the attacker knows the control input and the state, that is,
  $\I_k=\{u_1^k, x_1^k\}$.
\item the attacker knows the control input and delayed
  measurements from the sensor, that is,
  $\I_k=\{u_1^k, \y_1^{k-d}\}$ for some $d \geq 1$.
\end{enumerate}
\end{ex}

\textbf{Stealthiness of an attacker:} The attacker is constrained in
the input $\tilde{u}_{1}^{\infty}$ it replaces since it seeks to be
stealthy or undetected by the controller. If the controller is aware
that an attacker has replaced the correct control sequence
$u_{1}^{\infty}$ by a different sequence $\u_{1}^{\infty}$, it can
presumably switch to a safer mode of operation. Notions of
stealthiness have been proposed in the literature before. As an
example, for noiseless systems,~\cite{FP-FD-FB:10y} showed that
stealthiness of an attacker is equivalent to the existence of zero
dynamics for the system driven by the attack. Similar
to~\cite{FP-FD-FB:10y}, we seek to define the notion of stealthiness
without placing any restrictions on the attacker or the controller
behavior. However, we need to define a similar notion for stochastic
systems when zero dynamics may not exist. To this end, we pose the
problem of detecting an attacker by the controller as a (sequential)
hypothesis testing problem. Specifically, the controller relies on the
received measurements to decide the following binary hypothesis
testing problem:
\begin{align*}
  H_0: &  \text{ No attack is in progress (the controller receives } y_1^k\text{);}\\
  H_1: &  \text{ Attack is in progress (the controller receives } \y_1^k\text{).}
\end{align*}
For a given detector employed at the controller to select one of the
two hypotheses, denote the probability of false alarm (i.e., the
probability of deciding $H_1$ when $H_0$ is true) at time $k$ by
$p^{F}_k$, and the probability of correct detection (i.e., the
probability of deciding $H_1$ when $H_1$ is true) at time $k+1$ by
$p^D_k$. 
	
One may envisage that stealthiness of an attacker implies
$p^{D}_{k}=0$. However, as is standard in detection theory, we need to
consider both the quantities $p^F_k$ and $p^D_k$
simultaneously.\footnote{For instance, a detector that always declares
  $H_{1}$ to be true will achieve $p_{k}^{D}=1$. However, it will not
  be a good detector because $p_{k}^{F}=1$.} In fact, intuitively, an
attack is harder to detect if the performance of \emph{any} detector
is independent of the received measurements. In other words, we define
an attacker to be stealthy if there exists no detector that can
perform better (in the sense of simultaneously achieving higher
$p_{k}^{D}$ and lower $p_{k}^{F}$) than a detector that makes a
decision by ignoring all the measurements and making a random guess to
decide between the hypotheses. We formalize this intuition in the
following definition.

\begin{defi}{\bf \emph{(Stealthy attacks)}}\label{defi:strict_stealthiness}
  Consider the problem formulation stated in
  Section~\ref{sec:sys_model}.  An attack $\u_1^\infty$ is
  \begin{enumerate}
  \item strictly stealthy, if there exists no detector such that
    $p^F_k < p^D_k$ for any $k>0$.
		%
  \item $\epsilon$-stealthy, if, given $\epsilon>0$
      and $0<\delta<1$, for any detector for which
      $0 < 1 - p^D_k \leq \delta$ for all times~$k$, it holds that
    \begin{equation}\label{eq:eps_stealthy_exponent}
      \limsup_{k\to\infty} -\frac{1}{k} \log p^F_k \leq \epsilon.
    \end{equation}

  \end{enumerate}
\end{defi}
Intuitively, an attack is strictly stealthy if no detector can perform
better than a random guess in deciding whether an attack is in
progress. Further, an attack is $\epsilon$-stealthy if there exists no
detector such that $0 < 1 - p^D_k \leq \delta$ for all time $k$ and
$p^F_k$ converges to zero exponentially fast with rate greater than
$\epsilon$ as $k\to\infty$.

\textbf{Performance metric:} The requirement to stay stealthy clearly
curtails the performance degradation that an attacker can cause. The
central problem that we consider is to characterize the worst
performance degradation that an attacker can achieve for a specified
level of stealthiness. 
In the presence of an attack (and if the controller is unaware of the
attack), it uses the corrupted measurements $\y_1^\infty$ in the
Kalman filter. Let $\hat{\tilde{x}}_1^\infty$ be the estimate of the
Kalman filter \eqref{eq:kf} in the presence of the attack
$\u_1^\infty$, which is obtained from the recursion
\begin{align*}
  \hat{\tilde{x}}_{k+1} = A \hat{\tilde{x}}_k + K \tilde{z}_k + Bu_k ,
\end{align*}
where the innovation is
$\tilde{z}_{k}\triangleq \tilde{y}_{k}-C\hat{\tilde{x}}_{k}$.
Note that the estimate $\hat{\x}_{k+1}$ is a sub-optimal MMSE estimate
of the state $x_k$ since it is obtained by assuming the nominal
control input $u_{k}$, whereas the system is driven by the attack
input $\u_{k}$. Also, note that the random sequence $\z_1^\infty$ need
neither be zero mean, nor white or Gaussian. 

Since the Kalman filter estimate depends on the measurement sequence
received, as a performance metric, we consider the covariance of the
error in the predicted measurement $\hat{\tilde y}_k$ as compared to
true value $y_k$. Further, to normalize the relative impact of the
degradation induced by the attacker among different components of this
error vector, we weight each component of the error vector by an
amount corresponding to how accurate the estimate of this component
was without attacks. Thus, we consider the performance index
\begin{align*}
  \E \left[ \left( \hat{\tilde y}_k - y_k \right)^T \Sigma_z^{-1}
  \left( \hat{\tilde y}_k - y_k \right) \right] = \text{Tr} (\tilde P_k W),
\end{align*}
where $\tilde P_k$ is the error covariance matrix in the presence of
an attack,
$\tilde{P}_{k} = \mathbb{E}\big[(\hat{\x}_{k} - x_{k})(\hat{\x}_{k} -
x_{k})^T\big]$, and $W = C^T \Sigma_z^{-1} C$. To
obtain a metric independent of the time and focus on the long term
effect of the attack, we consider the limit superior of the arithmetic
mean of the sequence $\{\tr(\P_k W)\}_{k=1}^\infty$ and define
\begin{align}\label{eq:tilde_P_def}
  \tilde{P}_W &\triangleq \limsup_{k\to\infty} \frac{1}{k} \sum_{n=1}^k \tr( \tilde{P}_{n}W).
\end{align}
Notice that if the sequence $\{\tr(\tilde{P}_k W)\}_{k=1}^\infty$ is
convergent, then $\lim_{k\to\infty} \tr(\tilde{P}_k W) = \tilde{P}_W$, which
equals the Ces\`{a}ro mean of
$\tilde{P}_{k}W$. 
	

\textbf{Problems considered in the paper:} We assume
that the attacker is interested in staying stealthy or undetected
for as long as possible while maximizing the error covariance
$\tilde{P}_{W}$. We consider two problems:
\begin{enumerate}
\item What is a suitable metric for stealthiness of an attacker in
  stochastic systems where Assumption~\ref{asmp_zero} holds? We
  consider this problem in Section~\ref{sec:stealthy}.
\item For a specified level of stealthiness, what is the worst
  performance degradation that an attacker can achieve? We consider
  this problem in Section \ref{sec:fundamental_limits_given_W}.
\end{enumerate}

\section{Stealthiness in Stochastic systems}\label{sec:stealthy}
Our first result provides conditions that can be used to verify if an
attack is stealthy or not.
\begin{thm}{\bf \emph{(KLD and stealthy attacks)}}\label{thm: KLD and stealthy}
Consider the problem formulation stated in Section~\ref{sec:sys_model}. An attack $\u_1^\infty$ is
\begin{enumerate}
\item strictly stealthy if and only if  $D\big(\y_1^k\big\|y_1^k\big)=0$ for all $k>0.$
\item $\epsilon$-stealthy if the corresponding observation sequence $\y_1^\infty$ is ergodic and satisfies
  \begin{equation}\label{eq:Chernoff-Stein}
    \lim_{k\to\infty} \frac{1}{k} D\big(\y_1^k \big\| y_1^k \big) \leq \epsilon.
  \end{equation}
\item  $\epsilon$-stealthy only if the corresponding observation sequence $\y_1^\infty$ satisfies~(\ref{eq:Chernoff-Stein}).
\end{enumerate}
\end{thm}

\begin{pf}
Presented in Appendix~\ref{apx:lem:Chernoff_Stein_conv}.
\qed\end{pf}


The following result provides a characterization of
$D\big(\y_1^k \big\| y_1^k \big)$ that contains additional insight
into the meaning of stealthiness of an attacker.
\begin{prop}{\bf \emph{(KLD and differential
      entropy)}}\label{prop:stealthiness_interpretation}
  The quantity $D\big(\y_1^k \big\| y_1^k \big)$ can be calculated as
\begin{equation}
 \frac{1}{k} D\big(\y_1^k \big\| y_1^k \big)= \frac{1}{k}\sum_{n=1}^k \Big( I\big(\z_1^{n-1} ; \z_n \big) + D\big(\z_n \big\| z_n\big) \Big),
\end{equation}
where $I\big(\z_1^{n-1}; \z_n \big)$ denotes the mutual information
between $\z_1^{n-1}$ and $\z_n$~\cite[Section 8.5]{Cover06}.
\end{prop}
\begin{pf}
  Due to the invariance property of the Kullback-Leibler divergence
  \cite{Kullback97}, we have
  \begin{align*}
    D\big(\y_1^k \big\| y_1^k \big) =  D\big(\z_1^k \big\| z_1^k \big) ,
  \end{align*}
  for every $k>0$. Further, note that $z_1^\infty$ is an i.i.d. sequence
  of Gaussian random vectors with $z_k\sim\N(0,\Sigma_z)$. From
  (\ref{eq:KLD}), we obtain
  \begin{align*}
\frac{1}{k}  D\big(\z_1^k \big\| z_1^k \big)
    %
      &\stackrel{(a)}{=} -\frac{1}{k} h\big(\z_1^k\big) - \frac{1}{k}\sum_{n=1}^k\E\big[ \log f_{z_n}(z_n) \big] \\
                                    & \stackrel{(b)}{=}  \frac{1}{k}\sum_{n=1}^k \Big( -h\big(\z_n\big|\z_1^{n-1}\big) + h(\z_n)\notag\\
                                    &\qquad\qquad - h(\z_n) - \E\big[ \log f_{z_n}(z_n) \big]\Big)\\
                                    & = \frac{1}{k}\sum_{n=1}^k \Big( I\big(\z_1^{n-1} ; \z_n \big) + D\big(\z_n \big\| z_n\big) \Big), 
  \end{align*}
  where $h\big(\z_1^k\big)$ is the differential entropy of $\z_1^k$,
  $I\big(\z_1^{n-1}; \z_n \big)$ denotes the mutual information
  between $\z_1^{n-1}$ and $\z_n$. Equality (a) holds because
  $z_1^\infty$ is an independent random sequence, while equality (b)
  follows by applying the chain rule of differential entropy
  \cite[Theorem~8.6.2]{Cover06} on the term
  $-\frac{1}{k}h\left(\tilde{z}_{1}^{k}\right)$ to obtain
  $\frac{1}{k}\sum_{n=1}^{k}-h\left(\tilde{z}_{n}|\tilde{z}_{1}^{n-1}\right)$,
  and adding and subtracting $h(\tilde{z}_{n})$.
\qed\end{pf}
Intuitively, the mutual information $I\big(\z_1^{n-1}; \z_n \big)$
measures how much information about $\z_n$ can be obtained from
$\z_1^{n-1}$, that is, it characterizes the memory of the sequence
$\z_1^\infty$. Similarly, the Kullback-Leibler divergence
$D(\z_n \| z_n)$ measures the dissimilarity between the marginal
distributions of $\z_n$ and
$z_n$. Proposition~\ref{prop:stealthiness_interpretation} thus states
that the stealthiness level of an ergodic attacker can be degraded in
two ways: (i) if the sequence $\z_1^\infty$ becomes autocorrelated, and
(ii) if the marginal distributions of the random variables $\z(k)$ in
the sequence $\z_1^\infty$ deviate from $\N(0, \Sigma_z)$.

\section{Fundamental Performance
  Limitations}\label{sec:fundamental_limits_given_W}
We are interested in the maximal performance degradation $\P_W$ that
an $\epsilon$-stealthy attacker may induce. We begin by proving a
converse statement that gives an upper bound for $\P_W$ induced by an
$\epsilon$-stealthy attacker in Section \ref{sec: converse}. In
Section \ref{sec: achievability right} we prove a tight achievability
result that provides an attack that achieves the upper bound when the
system $(A,B,C)$ is right-invertible. In Section
\ref{subsec:achievability_nri} we prove a looser achievability result
that gives a lower bound on the performance degradation for non
right-invertible systems.

\subsection{Preliminary results}\label{sec: converse}
We will use a series of preliminary technical results, which will be
used in the following sections to present the main results of the
paper.

\begin{lem}
\label{lem:vijay_2}
Define the function $\bar\delta: [0,\infty)\to [1, \infty)$  as
  \begin{align}\label{eq:delta_bar}
    \bar{\delta}(x) = 2x +1 + \log \bar{\delta}(x).
  \end{align}
Then, for any $\gamma>0$, 
\begin{align}\label{eq:log_det_opt}
\bar{\delta}(\gamma) =&\arg\max_{x\in\R}  \; x,  \\
&\text{subject to } \frac{1}{2}x-\gamma-\frac{1}{2}\leq \frac{1}{2}\log x
  \notag .
\end{align}
\end{lem}
\begin{pf} 
Since a logarithm function is concave, the feasible
region of $x$ in (\ref{eq:log_det_opt}) is a closed interval upper
bounded by $\bar{\delta}(\gamma)$ as defined in \eqref{eq:delta_bar}. 
Thus, the result follows.
\qed\end{pf}
The following result is proved in the appendix.
\begin{lem}\label{lem:vijay_1}
  Consider the problem setup above. We have
  \begin{multline}
    \frac{1}{2k}\sum_{n=1}^k \tr\big(
    \E\big[\z_n\z_n^T\big]\Sigma_z^{-1} \big) \leq \frac{N_y}{2}
    +\frac{1}{k} D\big(\z_1^k \big\| z_1^k \big) + \\ \frac{N_y}{2}
    \log \Bigg( \frac{1}{N_y k}\sum_{n=1}^k
    \tr\big(\E[\z_n\z_n^T]\Sigma_z^{-1}\big)\Bigg) \label{eq:AM_GM_MIMO}.
  \end{multline}
  Further, if the sequence $\z_1^\infty$ is a sequence of independent
  and identically distributed (i.i.d.)  Gaussian random variables,
  $\z_{k}$, each with mean zero and covariance matrix
  $\E\big[\z_k\z_k^T\big]=\alpha \Sigma_z$, for some scalar $\alpha$,
  then~(\ref{eq:AM_GM_MIMO}) is satisfied with equality.
\end{lem}
\begin{pf}
See Appendix~\ref{apx:vijay_1}.
\qed\end{pf}
Combining Lemmas~\ref{lem:vijay_2} and~\ref{lem:vijay_1} leads to the following result.
\begin{lem}
\label{lem:vijay_3}
Consider the problem setup above. We have
 \begin{equation}
    \frac{1}{N_y k}\sum_{n=1}^k \tr\big( \E\big[\z_n\z_n^T\big]\Sigma_z^{-1} \big) \leq \bar{\delta}\Big(\frac{1}{N_y k}   D\big(\z_1^k \big\| z_1^k \big) \Big),
    \end{equation}
    where $\bar{\delta}(.)$ is as defined in \eqref{eq:delta_bar}.
    \end{lem}
    \begin{pf}
Proof follows from Lemma~\ref{lem:vijay_2} by using~(\ref{eq:AM_GM_MIMO}) and substituting
\begin{align*}
x&=\frac{1}{N_y k}\sum_{n=1}^k \tr\big( \E\big[\z_n\z_n^T\big]\Sigma_z^{-1} \big)\\
\gamma&=\frac{1}{N_y k}   D\big(\z_1^k \big\| z_1^k \big).
\end{align*}
\qed\end{pf}
The following result relates the covariance of the innovation and the
observation sequence.
\begin{lem}
Consider the problem setup above. We have
\begin{align}
\label{eq:real_inn}
CP_{k}C^{T}&=\E\big[z_{k}z_{k}^{T}\big]-\Sigma_{v}\\
\label{eq:attack_inn}
C\tilde{P}_{k}C^{T}&=\E\big[\z_{k}\z_{k}^{T}\big]-\Sigma_{v}.
\end{align}
\end{lem}
\begin{pf}
By definition,
\begin{equation}
z_k = y_k - C\hat{x}_k = C(x_k - \hat{x}_k) + v_k, 
\end{equation}
and similarly
\begin{equation}
\z_k=C(\tilde x_k - \hat{\x}_k) + v_k.
\end{equation}
Now both $(x_k - \hat{x}_k)$ and $(\tilde x_k - \hat{\x}_k)$ are
independent of the measurement noise $v_k$ due to Assumptions 1 and
7. Thus, we have
  \begin{align*}
CP_{k}C^{T}+\Sigma_{v}&=\E\big[z_{k}z_{k}^{T}\big]\\
C\tilde{P}_{k}C^{T}+\Sigma_{v}&=\E\big[\z_{k}\z_{k}^{T}\big],
\end{align*}
and the result follows.
\qed\end{pf} 

\subsection{Converse}\label{sec: converse}
We now present an upper bound of the weighted MSE induced by an
$\epsilon$-stealthy attack.
\begin{thm}[{\bf\emph{Converse}}]\label{thm:converse_MIMO}
  Consider the problem setup above. For any $\epsilon$-stealthy attack
  $\u_1^\infty$ generated by an information pattern $\I_1^\infty$ that
  satisfies Assumptions~\ref{assm_A1}-\ref{assm-A3}, we have
  \begin{align}\label{eq:converse_MIMO}
    \P_W  &\leq \tr (PW) + \left(
            \bar{\delta} \!\left(\frac{\epsilon}{N_y}\right)-1 \right) N_y,
  \end{align}
  where $N_y$ is the number of outputs of the system, the function
  $\bar\delta$ is defined in~(\ref{eq:delta_bar}), and $\tr(PW)$ is
  the weighted MSE in the absence of the attacker.
\end{thm}
\begin{pf}
We begin by writing
\begin{align*}
\P_{W}&=\limsup_{k\rightarrow\infty}\frac{1}{k}\sum_{n=1}^{k}\tr(\P_{n} C^{T}\Sigma_{z}^{-1}C)\\
&=\limsup_{k\rightarrow\infty}\frac{1}{k}\sum_{n=1}^{k}\tr(C\P_{n} C^{T}\Sigma_{z}^{-1})\\
&= \limsup_{k\to\infty} \frac{1}{k} \sum_{n=1}^k \tr\Big(\big(\E[\z_n \z_n^T] - \Sigma_v \big)\Sigma_z^{-1}\Big),
\end{align*}
where we have used the invariance of trace operator under cyclic
permutations and the relation in~(\ref{eq:attack_inn}),
respectively. The right hand side has two terms. The first term can be
upper bounded using Lemma~\ref{lem:vijay_3}, so that we obtain
\begin{equation*}
  \P_{W}\leq\limsup_{k\to\infty} N_y\bar{\delta} \! \left(\frac{1}{N_y k} D\!\left(\z_1^k \big\| z_1^k \right)\right) - \tr\big(\Sigma_v\Sigma_z^{-1}\big).
\end{equation*}
Since the function $\bar{\delta}$ is continuous and monotonic, we can
rewrite the above bound as
\begin{equation*}
  \P_{W}\leq N_y\bar{\delta} \! \left(\limsup_{k\to\infty}\frac{1}{N_y k} D\big(\z_1^k \big\| z_1^k \big)\right) - \tr\big(\Sigma_v\Sigma_z^{-1}\big).
\end{equation*}
Since the attack is $\epsilon$-stealthy, we use
Theorem~\ref{thm: KLD and stealthy} to bound the Kullback-Leibler
divergence $D\big(\z_1^k \big\| z_1^k \big)$ to obtain
\begin{equation}
\P_{W}\leq  N_y\bar{\delta}\Big(\frac{\epsilon}{N_y}\Big) -
  \tr\big(\Sigma_v\Sigma_z^{-1}\big).
  \end{equation}
  Finally, substituting for $\Sigma_{v}$ from~(\ref{eq:real_inn}) on
  the right hand side and using $W = C^T\Sigma_z^{-1} C$ completes the
  proof.
  \qed
\end{pf}

\begin{rem}{\bf \emph{(Stealthiness vs induced error)}}
  Theorem~\ref{thm:converse_MIMO} provides an upper bound for the
  performance degradation $\P_W$ for $\epsilon$-stealthy
  attacks. Since $\bar\delta \!\left(\frac{\epsilon}{N_y} \right)$ is
  a monotonically increasing function of $\epsilon$, the upper bound
  (\ref{eq:converse_MIMO}) characterizes a trade-off between the
  induced error and the stealthiness level of an attack. 
\end{rem}

To further understand this result, we consider two extreme cases,
namely, $\epsilon=0$, which implies strictly stealthiness, and
$\epsilon\to\infty$, that is, no stealthiness level.
\begin{cor}
A strictly stealthy attacker cannot induce any performance degradation.
\end{cor}
\begin{pf}
  A strictly stealthy attacker corresponds to $\epsilon=0$. Using the
  fact that $\bar{\delta}(0)=1$ in Theorem~\ref{thm:converse_MIMO}
  yields that $\tr(\P W)\leq \tr(PW)$.
\qed\end{pf}

\begin{cor}\label{cor:eps_infty}
  Let the attacker be $\epsilon$-stealthy. The upper bound in
  (\ref{eq:converse_MIMO}) increases linearly with $\epsilon$ as
  $\epsilon\to\infty$.
\end{cor}
\begin{pf}
  The proof follows from Theorem~\ref{thm:converse_MIMO} by noting
  that the first order derivative of the function $\bar\delta(x)$ is
  given by
  \begin{align}\label{eq:d_delta}
    \frac{d \bar\delta(x)}{d x}  &= \frac{d}{d x} \Big( 2x + 1 + \log \bar\delta(x) \Big) \notag\\
                                 &= 2 + \frac{1}{\bar\delta(x)} \frac{d \bar\delta(x)}{d x}\notag\\	
                                 & =
                                   \frac{2}{1-\frac{1}{\bar\delta(x)}}
                                   \rightarrow 2,
  \end{align}
  from the right as $x$ tends to infinity. 
\qed\end{pf}

\subsection{Achievability for Right Invertible Systems}\label{sec:
  achievability right}
 In this section, we show that the bound presented in Theorem~\ref{thm:converse_MIMO} is achievable if the system $(A,B,C)$ is right invertible. We begin with the following preliminary result.
 \begin{lem}
Let the system $(A,B,C)$ be right invertible. Then, the system $(A-KC, B, C)$ is also right  invertible.  
\end{lem}
\begin{pf}
  From Assumption 2, the system $(A,B,C)$ has no invariant zero. Since
  $(A,B,C)$ is also right invertible, the dimension $N_{u}$ of the
  control vector $u_{k}$ is no less the dimension $N_{y}$ of the
  output vector $y_{k}$. Now $(A-KC, B, C)$ is generated from
  $(A,B,C)$ using output feedback, and hence the system $(A-KC, B, C)$
  does not have invariant zeros either. Since the dimension of
  the input and output vectors of this new system remains $N_{u}$ and
  $N_{y}$ with $N_u\geq N_y$, the system $(A-KC, B,C)$ is right
  invertible.
\qed\end{pf}

Let $G_{RI}'$ be the right inverse of the system $(A-KC, B, C)$. Specifically, we will consider the following attack. 
  
\textbf{Attack $\mathcal{A}_{1}$:} The attack sequence is generated in
three steps. In the first step, a sequence $\zeta_1^\infty$ is
generated, such that each vector $\zeta_{k}$ is independent and
identically distributed and independent of the information pattern
$\I_{k}$ of the attacker, with probability density function
$\zeta_k\sim\N\big(0, \big( \bar{\delta}(\frac{\epsilon}{N_y}) -1
\big)\Sigma_z \big).$
In the second step, the sequence $\phi_1^\infty$ is generated as the
output of the system $G_{RI}'$ with $\zeta_1^\infty$ as the input
sequence. Finally, the attack sequence $\u_1^\infty$ is generated as
  \begin{equation}
  \label{eq:opt_attack_RI}
    \u_k = u_k + \phi_k.
  \end{equation}	

  \begin{rem}{\bf \emph{(Information pattern of attack $\mathcal{A}_{1}$)}}
    Notice that the attack $\mathcal{A}_{1}$ can be generated by an
    attacker with any information pattern satisfying Assumptions
    5--7.
  \end{rem}

We note the following property of the attack $\mathcal{A}_{1}$.
\begin{lem}\label{lem:vijay_4}
  Consider the attack $\mathcal{A}_{1}$. With this attack, the
  innovation sequence $\z_1^\infty$ as calculated at the controller,
  is a sequence of independent and identically distributed Gaussian
  random vectors with mean zero and covariance matrix
  \begin{equation}
    \label{eq:tilde_z_cov}
    \E[\z_k\z_k^T]  
    = \bar{\delta} \!\left( \frac{\epsilon}{N_y} \right)\Sigma_z.
  \end{equation}
\end{lem}
\begin{pf}
 %
  Consider an auxiliary Kalman filter that is implemented as the
  recursion
  \begin{align}\label{eq: za}
    \hat{x}_{k+1}^a = A \hat{x}_k^a + Kz_k^a+ B\u_k,
  \end{align}
  with the initial condition $\hat{x}_1^a=0$ and the innovation
  $z_k^a = \y_k-C\hat{x}_k^a$. The innovation sequence is independent
  and identically distributed with each $z_k^a\sim\N(0,
  \Sigma_z)$.  
  Now, we express $\z_k$ as
  \begin{align}
    \nonumber		\z_k &= \y_k - C\hat{\x}_k \\
    \nonumber		&= \y_k - C\hat{x}_k^a +  C(\hat{x}_k^a - \hat{\x}_k) \\
    \label{eq:z_k}		&= z_k^a - C\e_{k}  ,
  \end{align}
  where $\e_{k}\triangleq\hat{\x}_k - \hat{x}_k^a$.  Further,
  $\e_{k}$ evolves according to the recursion
  \begin{align}
    \nonumber		\e_{k+1} &= (A\hat{\x}_k + K \z_k + Bu_k) - (A\hat{x}_k^a + K z_k^a + B\u_k)\\
    \nonumber &=(A-KC) \e_k + B(u_k - \u_k) \\
    \label{eq:e_k}		& = (A-KC) \e_k - B\phi_{k} ,
  \end{align}
  with the initial condition $\e_1=0$. Together,~(\ref{eq:z_k})
  and~(\ref{eq:e_k}) define a system of the form
  \begin{equation}\label{eq:error_dyn_MIMO}
    \begin{aligned}
      \e_{k+1} &= (A-KC) \e_k + B(-\phi_{k}) ,\\
      z_{k}^{a}-\z_k &= C \e_k.
    \end{aligned}		
  \end{equation}
  We now note that (i) the above system is $(A-KC,B,C)$, (ii)
  $\phi_1^\infty$ is the output of the right inverse system of
  $(A-KC, B, C)$ with input $\zeta_1^\infty$, and (iii) the system in
  equation~(\ref{eq:error_dyn_MIMO}) is linear. These three facts
  together imply that the output of~(\ref{eq:error_dyn_MIMO}), i.e.,
  $\{z_{k}^{a}-\z_k\}_{k=1}^{\infty} $ is a sequence of independent
  and identically distributed random variables with each random
  variable distributed as
  $\N\big(0, \big( \bar{\delta} (\frac{\epsilon}{N_y}) -1 \big)\Sigma_z
  \big)$.
  Now since $z_k^a$ is independent of $\e_1^{k}$, we obtain that
  $\z_1^\infty$ is an independent and identically distributed sequence
  with each random variable $\z_{k}$ as Gaussian with mean zero and
  covariance matrix
  \begin{align*}
    \E[\z_k\z_k^T] = \left(\bar{\delta} \!\left( \frac{\epsilon}{N_y} \right) -1
    \right)\Sigma_z + \Sigma_z =
    \bar{\delta} \!\left(\frac{\epsilon}{N_y} \right)\Sigma_z. 
  \end{align*}
\qed\end{pf}
We now show that the attack $\mathcal{A}_{1}$ achieves the converse result in Theorem~\ref{thm:converse_MIMO}.
\begin{thm}{\bf\emph{(Achievability for right invertible
      systems)}}\label{thm:achievability}
  Suppose that the LTI system $(A,B,C)$ is right invertible. The
  attack $\mathcal{A}_{1}$ is $\epsilon$-stealthy and, with this attack,
  \begin{equation*}
    \P_W 
    = \tr(PW) + N_y\left(\bar{\delta} \! \left(\frac{\epsilon}{N_y}\right)-1 \right),
  \end{equation*}
  where $W=C^T\Sigma_z^{-1}C$. 
\end{thm}
\begin{pf}
  For the attack $\mathcal{A}_{1}$, Lemma~\ref{lem:vijay_4} states
  that $\z_1^\infty$ is a sequence of independent and identically
  distributed (i.i.d.) Gaussian random variables $\z_{k}$ each with
  mean zero and covariance matrix
  $\E\big[\z_k\z_k^T\big]=\alpha \Sigma_z,$ with
  $\alpha=\bar{\delta}(\frac{\epsilon}{N_y}).$
  Lemma~\ref{lem:vijay_1}, thus, implies that~(\ref{eq:AM_GM_MIMO})
  holds with equality.
  Further, following the proof of Theorem~\ref{thm:converse_MIMO},
  if~(\ref{eq:AM_GM_MIMO}) holds with equality,
  then~(\ref{eq:converse_MIMO}) also holds with equality. Thus, the
  attack $\mathcal{A}_{1}$ achieves the converse in terms of
  performance degradation.
  	
  Next we show that the attack is $\epsilon$-stealthy. Once again,
  from Lemma~\ref{lem:vijay_1} and~(\ref{eq:tilde_z_cov}), we have for
  every $k>0$,
  \begin{align*}
    \frac{1}{k} D\big( \z_1^k\|z_1^k \big) &= 
                                             \frac{1}{2k}\sum_{n=1}^k \tr\big( \E\big[\z_n\z_n^T\big]\Sigma_z^{-1} \big) - \frac{N_y}{2} \\&\qquad- \frac{N_y}{2} \log  \Bigg( \frac{1}{N_y k}\sum_{n=1}^k \tr\big(\E[\z_n\z_n^T]\Sigma_z^{-1}\big)\Bigg) \\
                                           &=	\frac{1}{2k}\sum_{n=1}^k \tr\Big( \bar{\delta}(\frac{\epsilon}{N_y})\Sigma_z \Sigma_z^{-1} \Big) - \frac{N_y}{2} \\
                                           &\qquad\qquad\qquad - \frac{1}{2k} \sum_{n=1}^k \log \det\big(\bar{\delta}(\frac{\epsilon}{N_y})\Sigma_z\Sigma_z^{-1}\big)\\
                                           & =\frac{N_y}{2}  \bar{\delta}\Big(\frac{\epsilon}{N_y}\Big) - \frac{N_y}{2} - \frac{N_y}{2}\log \bar{\delta}\Big(\frac{\epsilon}{N_y}\Big)\\
                                           & = \epsilon.
  \end{align*}
  Now with this attack, $\z_1^\infty$ is an independent and
  identically distributed sequence and the measurement sequence
  $\y_1^\infty$ is ergodic. Thus, from Theorem~\ref{thm: KLD and
    stealthy}, the attack $\mathcal{A}_{1}$ is $\epsilon$-stealthy.
\qed\end{pf}
	
\begin{rem}{\bf \emph{(Attacker information
      pattern)}}\label{rem:info_pattern}
  Intuitively, we may expect that the more information about the state
  variables that an attacker has, larger the performance degradation
  it can induce. However, Theorem \ref{thm:converse_MIMO} and Theorem
  \ref{thm:achievability} imply that the only critical piece of
  information for the attacker to launch an optimal attack is the
  nominal control input $u_1^\infty$.
\end{rem}

\subsection{Achievability for System that are not
  Right Invertible}\label{subsec:achievability_nri}
If the system is not right invertible, the converse result in
Theorem~\ref{thm:converse_MIMO} may not be achieved. We now construct a heuristic attack $\mathcal{A}_{2}$ that allows
us to derive a lower bound for the performance degradation $\P_W$
induced by $\epsilon$-stealthy attacks against systems that are not
right invertible.

\textbf{Attack $\mathcal{A}_{2}$:} The attack sequence is generated
as 
\begin{equation}\label{eq:attack_achievability_non_RI}
  \u_k = u_k + L\e_k - \zeta_k ,
\end{equation}
where $\e_k = \hat{\x}_k - \hat{x}_k^a$ as in
\eqref{eq:error_dyn_MIMO}, and the sequence $\zeta_1^\infty$ is generated
such that each vector $\zeta_{k}$ is independent and identically
distributed with probability density function
$\zeta_k\sim\N\big(0, \Sigma_\zeta \big)$ and independent of the information pattern $\I_{k}$ of the
attacker. The feedback matrix $L$ and
the covariance matrix $\Sigma_\zeta$ are determined in three steps,
which are detailed next.

\emph{Step 1 (Limiting the memory of the innovation sequence
  $\z_1^\infty$):} Notice that, with the attack $\mathcal{A}_2$ and
the notation in \eqref{eq: za}, the dynamics of $\e_k$ and $\z_k$ are given by
\begin{align}\label{eq:error_dyn_MIMO_achievability}
  \begin{aligned}
    \e_{k+1} &= (A-KC-BL) \e_k + B\zeta_k \\
    \z_k &= C \e_k + z_k^a.
  \end{aligned}
\end{align}
The feedback matrix $L$ should be selected to eliminate the memory of
the innovation sequence computed at the controller. One way to achieve this aim is to set $A-KC-BL=0$. In other words, if $A-KC-BL=0$, then $\z_1^\infty$
is independent and identically distributed. It may not be possible to select $L$ to achieve this aim exactly. Thus, we propose the following heuristic. Note that if $A-KC-BL=0$, then the cost function
\begin{equation}\label{eq:cheap_control_LQG}
  \lim_{k\to\infty} \frac{1}{k} 
  \sum_{n=1}^k \tr\big( \E[\e_n\e_n^T]W \big),
\end{equation}
is minimized, with $W = C^T\Sigma_z^{-1} C$. Since
\begin{align*}
  \sum_{n=1}^k \tr\big( \E[\e_n\e_n^T]W \big) 
                                              &= \E \left[ \sum_{n=1}^k 
                                                \e_n^T W
                                                \e_n \right],
\end{align*}
selecting $L$ to satisfy the constraint $A-KC-BL=0$ is equivalent to
selecting $L$ to solve a cheap Linear Quadratic Gaussian (LQG) problem
\cite[Section VI]{hespanha2009linear}. Thus, heuristically, we select
the attack matrix $L$ as the solution to this cheap LQG problem and,
specifically, as
\begin{equation}
  L = \lim_{\eta\rightarrow0}(B^T T_{\eta} B+\eta I)^{-1} B^T T_{\eta} (A-KC) ,
\end{equation}
where $T_{\eta}$ is the solution to the discrete 
algebraic Riccati equation
\begin{multline*}
  T_{\eta} = \\(A-KC)^T \Big(T_{\eta} -T_{\eta} B(B^TT_{\eta} B+\eta I)^{-1} B^TT_{\eta} \Big) (A-KC) + W .
\end{multline*}


\emph{Step 2 (Selection of the covariance matrix $\Sigma_\zeta$):}
Notice that the selection of the feedback matrix $L$ in Step 1 is
independent of the covariance matrix $\Sigma_\zeta$. As the second step, we select the covariance matrix $\Sigma_\zeta$ such that
$C\Sigma_\e C^T$ is close to a scalar multiplication of $\Sigma_z$,
say $\alpha^2 \Sigma_z$. From \eqref{eq:error_dyn_MIMO_achievability},
notice that
\begin{align*}
  \lim_{k\to\infty} \E[\z_k\z_k^T] = C\Sigma_{\e}C^T + \Sigma_z ,
\end{align*}
where $\Sigma_{\e}\in\mathbb{S}_{+}^{N_x}$ is the positive
semi-definite solution to the equation
\begin{equation}\label{eq:lyapunov_error}
  \Sigma_{\e} = (A-KC-BL) \Sigma_{\e} (A-KC-BL)^T + B\Sigma_\zeta B^T.
\end{equation}
We derive an expression for $\Sigma_\zeta$ from
\eqref{eq:lyapunov_error} by using the pseudoinverse matrices of $B$
and $C$, i.e.,
\begin{align}\label{eq:Sigma_zeta_pseudo}
  &\Sigma_\zeta = \alpha^2B^\dagger\Big( C^\dagger\Sigma_z
                  (C^T)^\dagger +  \\
                &-(A-KC-BL) C^\dagger\Sigma_z (C^T)^\dagger(A-KC-BL)^T \Big)(B^T)^\dagger,\notag
\end{align}	
where $^\dagger$ denotes the pseudoinverse operation. It should be
noted that the right-hand side of \eqref{eq:Sigma_zeta_pseudo} may not
be positive semidefinite. Many choices are possible to construct a
positive semi-definite $\Sigma_{\zeta}$. We propose that if the the
right-hand side is indefinite, we set its negative eigenvalues to zero
without altering its eigenvectors. Note that, if $B$ and $C$ are both
invertible, then we could directly set
$C\Sigma_\e C^T = \alpha^2 \Sigma_z$.


\emph{Step 3 (Enforcing the stealthiness level):} The covariance
matrix $\Sigma_\zeta$ obtained in Step 2 depends on the parameter
$\alpha$. In this last step, we select $\alpha$ so as to make the
attack $\mathcal A_2$ $\epsilon$-stealthy. To this aim, we first
compute an explicit expression for the stealthiness level and the
error induced by $\mathcal{A}_2$.

For the entropy rate of $\z_1^\infty$, since $\z_1^\infty$ is
Gaussian, we obtain
\begin{align}
  \lim_{k\to\infty}\frac{1}{k} h\big(\z_1^k\big) &= \lim_{k\to\infty} h\big(\z_{k+1} \big\vert \z_1^{k}\big)\label{eq:entropy_rate_stationary}\\
                                                 & = \lim_{k\to\infty} \frac{1}{2}\log \Big( (2\pi e)^{N_y} \det\big(\E[(\z_{k+1}-g_k(\z_1^k))\notag\\
                                                 &\qquad\qquad\qquad\qquad(\z_{k+1}-g_k(\z_1^k))^T] \big) \Big ) \label{eq:Gaussian_Fano}\\
                                                 & =  \frac{1}{2}\log \big( (2\pi e)^{N_y} \det(CSC^T + \Sigma_z) \big) \label{eq:Gaussian_entropy_rate}
\end{align}
where $g_k(\z_1^k)$ is the minimum mean square estimate of $\e_{k+1}$
from $\z_1^k$, which can be obtained from Kalman filtering, and
$S\in\mathbb{S}_{+}^{N_y}$ is the positive semidefinite solution to
the following discrete algebraic Riccati equation 
\begin{align}
  S &= (A-KC-BL)\Big( S - SC^T(CSC^T+\Sigma_z)^{-1} C S\Big) \notag\\
    &\qquad \times(A-KC-BL)^T + B\Sigma_\zeta B^T.\label{eq:reccati_entropy}
\end{align}
Note that the equality \eqref{eq:entropy_rate_stationary} is due
to~\cite[Theorem~4.2.1]{Cover06}; equality \eqref{eq:Gaussian_Fano} is
a consequence of the maximum differential entropy
lemma~\cite[Section 2.2]{ElGamal11}; the positive semidefinite matrix $S$ that
solves \eqref{eq:reccati_entropy} represents the steady-state error
covariance matrix of the Kalman filter that estimates $\z_{k+1}$ from
$\z_1^k$. Thus, the level of stealthiness for the attack
$\mathcal{A}_2$ is
\begin{align}\label{eq:KLD_achievability}
  &\lim_{k\to\infty}\frac{1}{k} D\big(\z_1^k \big\| z_1^k \big) =
    \epsilon = -\frac{1}{2}\log ( (2\pi e)^{N_y} \det(CSC^T + \Sigma_z)\notag\\
  & \qquad + \frac{1}{2}\log\big( (2\pi)^{N_y} \det(\Sigma_z) \big) + \frac{1}{2} \tr\big( (C\Sigma_\e C ^T + \Sigma_z)\Sigma_z^{-1} \big)\notag\\
  &= -\frac{1}{2} \log\det(I + SW) +\frac{1}{2} \tr(\Sigma_\e W) + \frac{1}{2} N_y,
\end{align}
where $W = C^T\Sigma_z^{-1}C$. To conclude our design of the attack
$\mathcal{A}_2$, we use \eqref{eq:KLD_achievability} to solve for the
desired value of $\alpha$, and compute the error induced by the
$\epsilon$-stealthy attack $\mathcal{A}_2$ as \begin{align}
                                                \P_W &= \lim_{k\to\infty} \frac{1}{k}\sum_{n=1}^k \tr(\E[\z_n\z_n^T]\Sigma_z^{-1}) - \tr(\Sigma_v\Sigma_z^{-1})\notag\\
                                                     & = \tr(PW) +
                                                       \tr(\Sigma_{\e}W)
                                                       -
                                                       N_y \label{eq:tilde_P_NRI}
\end{align}
where $\Sigma_\e$ is the solution to the Lyapunov equation
\eqref{eq:lyapunov_error}.

	\section{Numerical Results}\label{sec:numerical}
\begin{figure}
	\centering
	\includegraphics[width=1\linewidth]{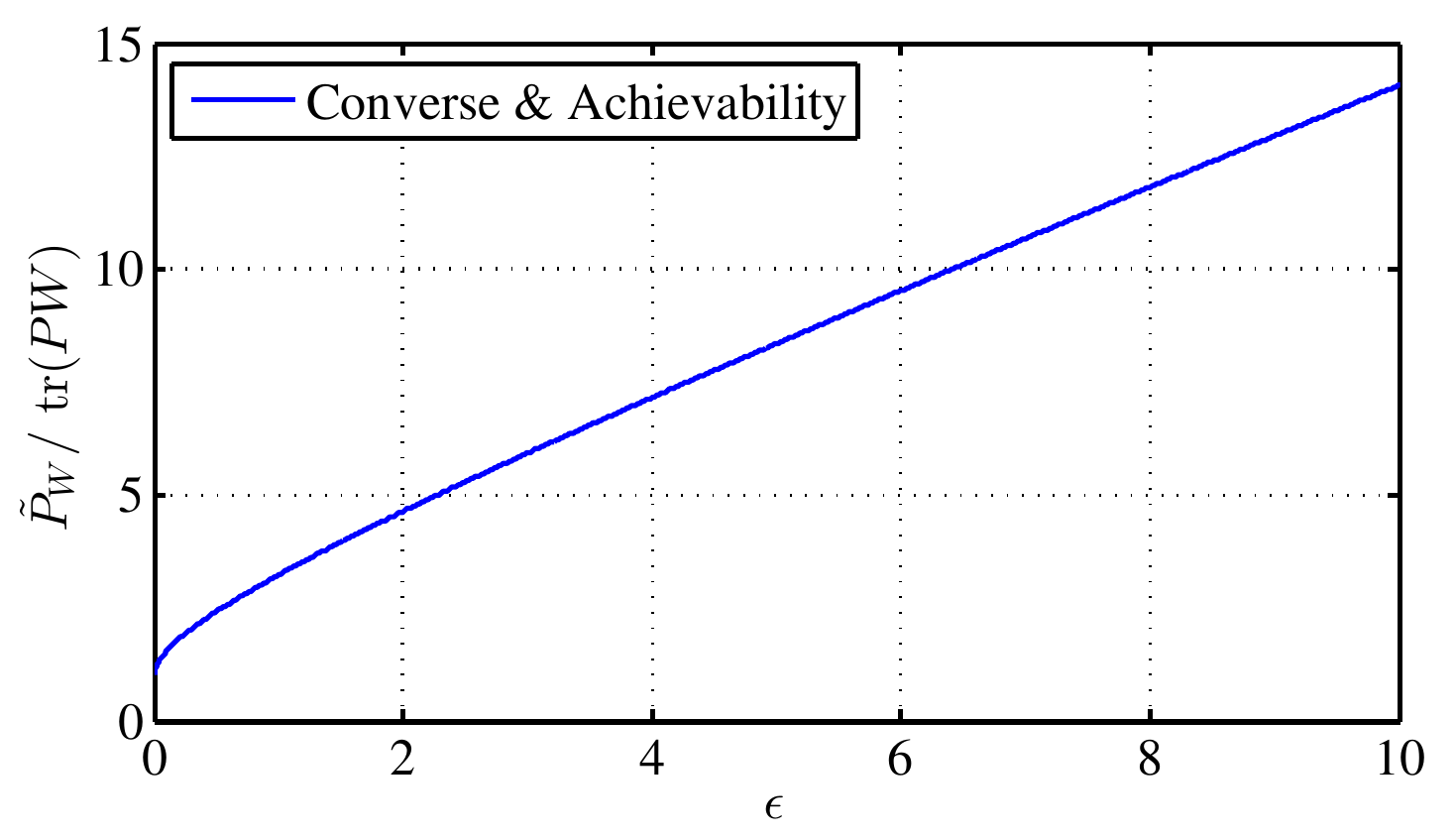}
	\caption{The converse and achievability for the right invertible system, where the weighted MSE $\P_W$ is the upper bound in (\ref{eq:converse_MIMO}) and the 	weight matrix $W=C^T\Sigma_z^{-1}C$.}
	\label{fig:fig_RI_P}
\end{figure}

\textbf{Example 1} Consider a right invertible system $(A,B,C)$ where
	\begin{align*}
	A &= \left[\begin{array}{cccc}
		2   &  0  &   0  &   0 \\
		0   &  -1  &   0  &   0 \\
		1   &  0  &  1  &   0 \\
		0   &  0  &   0  &   2 
	\end{array}\right], B=\left[
	\begin{array}{cc}
		1  &   0 \\
	    1  &   0 \\
     	0  &   2 \\
     	0  &   1 
     	\end{array} \right],   C=\left[\begin{array}{cc}
		0  &   0 \\
		0  &   1 \\
		2  &   0 \\
		0  &   1 
     	\end{array} \right]^T,
	\end{align*}
	and let $\Sigma_w = 0.5 I$ and $\Sigma_v = I$. Figure~\ref{fig:fig_RI_P} plots the upper bound (\ref{eq:converse_MIMO}) of performance degradation achievable for an attacker versus the attacker's stealthiness level $\epsilon$. From Theorem~\ref{thm:achievability}, the upper bound can be achieved by a suitably designed $\epsilon$-stealthy attack. Thus, Fig.~\ref{fig:fig_RI_P} represents a fundamental limitation for the performance degradation that can be induced by any $\epsilon$-stealthy attack. Observe that plot is approximately linear as $\epsilon$ becomes large, as predicted by Corollary~\ref{cor:eps_infty}.

\begin{figure}
	\centering
	\includegraphics[width=1\linewidth]{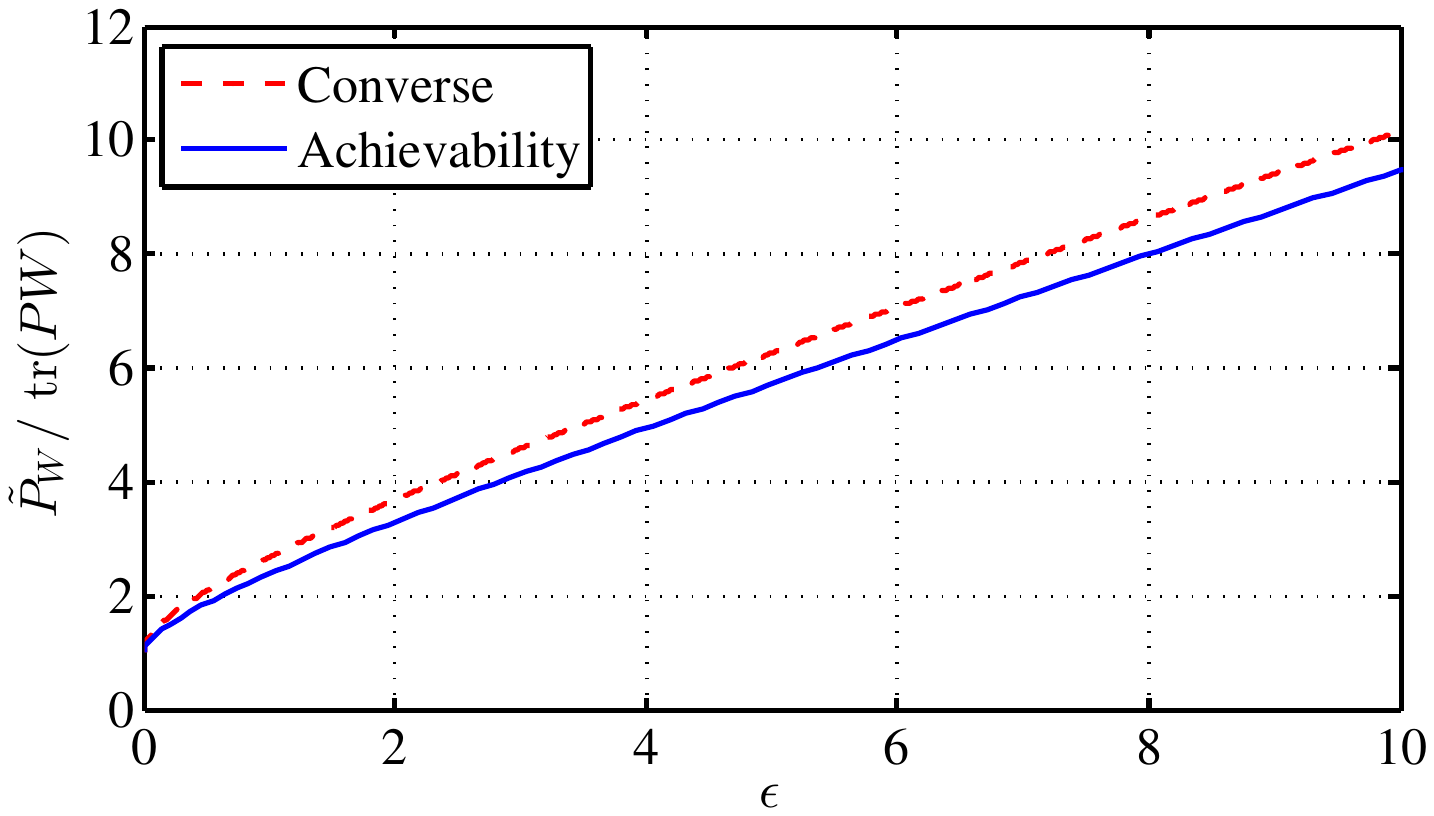}
	\caption{The converse and achievability for the right non-invertible system, where the weight matrix $W=C^T\Sigma_z^{-1}C$. The converse is obtained from (\ref{eq:converse_MIMO} and the achievability is the weighted MSE $\P_W$ induced by the heuristic algorithm $\mathcal{A}_{2}$.}
	\label{fig:fig_NRI_P}
\end{figure}

\textbf{Example 2}
Consider the system $(A,B,C)$ with
\begin{align*}
  A &= \left[\begin{array}{ccccc}
               2  &  -1  &   0  &   0  &   0 \\
               1  &  -3  &   0  &   0  &   0 \\
               0  &   0  &  -2  &   0  &   0 \\
               0  &   0  &   0  &  -1  &   0 \\
               0  &   0  &   0  &   0  &   3     
             \end{array}\right], \quad B=\left[
                                         \begin{array}{cc}
                                           2  &   0 \\
                                           1  &   0 \\
                                           0  &   1 \\
                                           0  &   1 \\
                                           1  &   1
                                         \end{array} \right], \\
  C &=\left[\begin{array}{ccccc}
              1  &  -1  &   2  &   0  &   0 \\
              -1  &    2  &   0  &   3  &   0 \\
              2   &  1  &   0  &   0 &    4
            \end{array} \right],
\end{align*}
which fails to be right invertible. Let $\Sigma_w = 0.5 I$ and
$\Sigma_v = I$. In Fig.~\ref{fig:fig_NRI_P}, we plot the upper bound
for the value of $\P_W$ that an $\epsilon$-stealthy attacker can
induce, as calculated using Theorem~\ref{thm:converse_MIMO}. The value
of $\P_W$ achieved by the heuristic attack $\mathcal{A}_{2}$ as
presented in Section~\ref{subsec:achievability_nri} is also
plotted. Although the bound is fairly tight as compared to the
performance degradation achieved by the heuristic attack; nonetheless,
there remains a gap between the two plots.

\section{Conclusion}\label{sec:conclusion}
This work characterizes fundamental limitations and achievability
results for performance degradation induced by an attacker in a
stochastic control system. The attacker is assumed to know the system
parameters and noise statistics, and is able to hijack and replace the
nominal control input. We propose a notion of $\epsilon$-stealthiness
to quantify the difficulty of detecting an attack from the
measurements, and we characterize the largest degradation of Kalman
filtering induced by an $\epsilon$-stealthy attack. For right
invertible systems, our study reveals that the nominal control input
is the only critical piece of information to induce the largest
performance degradation. For systems that are not right invertible, we
provide an achievability result that lower bounds the performance
degradation that an optimal $\epsilon$-stealthy attack can achieve.

\appendix
\section{Proof of Theorem~\ref{thm: KLD and
    stealthy}}\label{apx:lem:Chernoff_Stein_conv}
	
The first statement follows directly from the Neyman-Pearson
Lemma~\cite{Poor98}.

For the second statement, we apply the Chernoff-Stein Lemma for
ergodic measurements (see \cite{Polyanskiy12}) that states that for
any given attack sequence $\u_1^\infty$, for a given
$0<1-p^D_k\leq\delta$ where $0<\delta<1$, the best achievable decay
exponent of $p_k^F$ is given by
$\lim_{k\to\infty} \frac{1}{k} D\big(\y_1^k \big\| y_1^k \big)$. For
this attack sequence and with any detector, we obtain
\begin{equation*}
  \limsup_{k\to\infty} -\frac{1}{k}\log p_k^F \leq \lim_{k\to\infty} \frac{1}{k} D\big(\y_1^k \big\| y_1^k \big) \leq \epsilon.
\end{equation*}
Thus, by Definition~\ref{defi:strict_stealthiness}, the attack is
$\epsilon$-stealthy.

Finally, the proof for the third statement follows by
contradiction. Assume that (\ref{eq:Chernoff-Stein}) does not hold and
there exists an $\epsilon$-stealthy attack $\u_1^\infty$ such that
$\limsup_{k\to\infty} \frac{1}{k} D\big(\y_1^k \big\| y_1^k \big) >
\epsilon$.
Suppose that the detector employs the standard log-likelihood ratio
test with threshold $\lambda_k$ at every time $k+1$. Thus, the test is
given by
\begin{equation*}
  L_k(\eta_1^k)^{\stackrel{H_0}<}_{\stackrel{\geq}{H_1}} \lambda_k, \quad\text{where}\quad L_k(\eta_1^k) = \log \frac{f_{\y_1^k}(\eta_1^k)}{f_{y_1^k}(\eta_1^k)}
\end{equation*}
is the log-likelihood ratio and $\eta_1^k=y_1^k$
(resp. $\eta_1^k=\y_1^k$) if $H_0$ (resp. $H_1$) is true.  Define the
conditional cumulant generating function for the log-likelihood ratio
to be $g_{k|0}(s) = \log\E\big[ e^{sL_k} \big\vert H_0 \big]$ and
$g_{k|1}(s) = \log\E\big[ e^{sL_k} \big\vert H_1 \big]$. Note that
$g_{k|0}(s)=g_{k|1}(s-1)$. Let $\lambda_k$ be chosen to ensure that
$0<1-p^D_k\leq\delta$ for every $k>0$ (notice that such $\lambda_k$
always exists, because $p^D_k$ increases to one as $\lambda_k$
decreases to zero). Then, for any $s_k >0$, Chernoff's inequality
yields
\begin{equation*}
  p^F_k = \mathbb{P}[L_k \geq \lambda_k \vert H_0] \leq e^{-s_k \lambda_k + g_{k|0}(s_k)},
\end{equation*}
or
\begin{equation*}
  -\log p^F_k \geq s_k \lambda_k - g_{k|0}(s_k).
\end{equation*}
Further manipulation yields
\begin{align}
  -\log p^F_k	& \geq s_k \lambda_k - g_{k|1}(s_k-1) \notag\\
                & = s_k \lambda_k - \log \E\big[ e^{(s_k-1)L_k} \big\vert H_1 \big].
\end{align}
Now, by applying Jensen's inequality twice we obtain
\begin{align*}
  -\log p^F_k
  & \geq s_k \lambda_k + \log \E\big[ e^{-(s_k-1)L_k} \big\vert H_1 \big]\\
  & \geq s_k \lambda_k + \E[-(s_k-1)L_k | H_1 ]. 
\end{align*}
Finally, using 
$\E[L_k|H_1] = D\big(\y_1^k \big\| y_1^k \big)$ implies
\begin{equation}\label{eq:error_exp_bound}
  -\log p^F_k \geq D\big(\y_1^k \big\| y_1^k \big) + s_k\Big(\lambda_k - D\big(\y_1^k \big\| y_1^k \big)\Big).
\end{equation}	
Now, for any time index $k$ such that
$\frac{1}{k} D\big(\y_1^k \big\| y_1^k \big) > \epsilon$, let
\begin{equation}\label{eq:sk}
  s_k =  \frac{ D\big(\y_1^k \big\| y_1^k \big)- k\epsilon }{ 2\Big| D\big(\y_1^k \big\| y_1^k \big) - \lambda_k \Big|}.
\end{equation}
Using (\ref{eq:error_exp_bound}), (\ref{eq:sk}) and 
${\limsup_{k\to\infty} }\frac{1}{k} D\big(\y_1^k \big\| y_1^k \big) >
\epsilon$, we obtain
\begin{equation*}
  \limsup_{k\to\infty} -\frac{1}{k} \log p^F_k > \epsilon ,
\end{equation*}
which contradicts \eqref{eq:eps_stealthy_exponent}. Hence, the attack
cannot be stealthy, and the condition stated in
(\ref{eq:Chernoff-Stein}) must be true.

\section{Proof of Lemma~\ref{lem:vijay_1}}\label{apx:vijay_1}
By definition, we can write the Kullback-Leibler divergence can be written as
\begin{align*}
  D\big(\z_1^k \big\| z_1^k \big)&=\int_{-\infty}^\infty f_{\z_{1}^{k}}(t_{1}^{k})\log f_{\z_{1}^{k}}(t_{1}^{k}) dt_{1}^{k}\\&\qquad\qquad-\int_{-\infty}^\infty f_{\z_{1}^{k}}(t_{1}^{k})\log f_{z_{1}^{k}}(t_{1}^{k}) dt_{1}^{k}\\
                                 &=-h\big(\z_1^{k}\big)-\int_{-\infty}^\infty f_{\z_{1}^{k}}(t_{1}^{k})\log f_{z_{1}^{k}}(t_{1}^{k}) dt_{1}^{k}.
\end{align*}
Now, $z_{1}^{k}$ is the innovation sequence without any attack and is
thus an independent and identically distributed sequence of Gaussian
random variables with mean 0 and covariance $\Sigma_{z}$. Plugging
into the above equation yields
\begin{multline*}
D\big(\z_1^k \big\| z_1^k \big)=-h\big(\z_1^{k}\big)+\frac{k}{2}\log\Big( (2\pi)^{N_y} \det(\Sigma_z) \Big)\\+\frac{1}{2}\sum_{n=1}^k \tr\big( \E\big[\z_n\z_n^T\big]\Sigma_z^{-1} \big),
\end{multline*}
which we can rewrite as
\begin{multline}
  \frac{1}{2k}\sum_{n=1}^k \tr\big(
  \E\big[\z_n\z_n^T\big]\Sigma_z^{-1} \big) = \frac{1}{k} D\big(\z_1^k
  \big\| z_1^k \big)\\ - \frac{1}{2}\log\Big( (2\pi)^{N_y}
  \det(\Sigma_z) \Big) + \frac{1}{k} h\big(\z_1^{k}\big).
\end{multline}
We can upper-bound the right hand side by first using the
sub-additivity property of differential entropy \cite[Corollary
8.6.1]{Cover06}, and then further bounding the entropy $h(\z_{n})$
using the maximum differential entropy lemma \cite[Section
2.2]{ElGamal11} for multivariate random variables. Thus, we obtain
\begin{align*}
  &\frac{1}{2k}\sum_{n=1}^k \tr\big( \E\big[\z_n\z_n^T\big]\Sigma_z^{-1} \big)\notag\\
  &\quad \leq \frac{1}{k}  D\big(\z_1^k \big\| z_1^k \big) -   \frac{1}{2}\log\Big( (2\pi)^{N_y} \det(\Sigma_z) \Big) + \frac{1}{k} \sum_{n=1}^k h(\z_n)\\
  & \quad \leq \frac{1}{k}  D\big(\z_1^k \big\| z_1^k \big) -   \frac{1}{2}\log\Big( (2\pi)^{N_y} \det(\Sigma_z) \Big)\\
  &\qquad\qquad\qquad + \frac{1}{k} \sum_{n=1}^k \frac{1}{2} \log\big( (2\pi e)^{N_y} \det(\E[\z_n \z_n^T])\big),
\end{align*}           
with equality if the sequence $\z_1^k$ is an independent sequence of
random variables with each random variable $\z_n$ as Gaussian
distributed with mean zero for all $n$.  Straight-forward
algebraic manipulation on the last two terms yields
\begin{align*}
  &\frac{1}{2k}\sum_{n=1}^k \tr\big( \E\big[\z_n\z_n^T\big]\Sigma_z^{-1} \big)\notag\\
    & \quad \leq \frac{1}{k}  D\big(\z_1^k \big\| z_1^k \big) -   \frac{1}{2}\log\Big( (2\pi)^{N_y}\Big)  - \frac{1}{2}\log\Big( \det(\Sigma_z) \Big)\\
  &\qquad+ \frac{1}{k} \sum_{n=1}^k \frac{1}{2} \log\big( (2\pi e)^{N_y}\big) + \frac{1}{k} \sum_{n=1}^k \frac{1}{2} \log\big(  \det(\E[\z_n \z_n^T])\big) \\
  & \quad\leq \frac{1}{k}  D\big(\z_1^k \big\| z_1^k \big) -   \frac{1}{2}\log\Big( (2\pi)^{N_y}\Big) + \frac{1}{k} \sum_{n=1}^k \frac{1}{2} \log\big( (2\pi e)^{N_y}\big) \\
  &\qquad + \frac{1}{k} \sum_{n=1}^k \frac{1}{2} \log\big(  \det(\E[\z_n \z_n^T])\big)- \frac{1}{2}\log\Big( \det(\Sigma_z) \Big)\\
  &= \frac{1}{k}  D\big(\z_1^k \big\| z_1^k \big) +  \frac{N_{y}}{2} + \frac{1}{k} \sum_{n=1}^k \frac{1}{2} \log\big(  \det(\E[\z_n \z_n^T])(\det(\Sigma_z))^{-1} \Big).
\end{align*}      
We can further bound
\begin{align*}
  \nonumber  \det(\E[\z_n \z_n^T])(\det(\Sigma_z))^{-1}&=  \det(\E[\z_n \z_n^T]\Sigma_z^{-1})\\
                                                       &\leq \left(\frac{1}{N_{y}}\tr(\E[\z_n \z_n^T]\Sigma_z^{-1})\right)^{N_{y}},
\end{align*}
so that
\begin{multline}
  \frac{1}{2k}\sum_{n=1}^k \tr\big( \E\big[\z_n\z_n^T\big]\Sigma_z^{-1} \big)
   \leq \frac{1}{k}  D\big(\z_1^k \big\| z_1^k \big)  + \frac{N_y}{2}\\+ \frac{N_y}{2k} \sum_{n=1}^k \log \Big(\frac{1}{N_y}\tr\big(\E[\z_n\z_n^T]\Sigma_z^{-1}\big)\Big),\label{eq:tr_det_ineq_MIMO}
\end{multline}
with equality if the matrix $\E[\z_n \z_n^T]$ is a scalar
multiplication of $\Sigma_z$ for all $n$ Finally, using the Arithmetic
Mean and Geometric Mean (AM-GM) inequality yields the desired result
\eqref{eq:AM_GM_MIMO}. For the AM-GM inequality to hold with equality
we need that $\tr(\E[\z_n\z_n^T]\Sigma_z^{-1})$ is constant for every
$n$. Collecting all the above conditions for equality at various
steps, \eqref{eq:AM_GM_MIMO} holds with equality if
$\E\big[\z_k\z_k^T\big]=\alpha \Sigma_z$ for some scalar $\alpha$.

\bibliography{./2015_TAC_bib,./BIB}
	
\end{document}